\documentclass[11pt, final]{article}
 \usepackage{a4}
 \usepackage{amsmath}%
 \usepackage{amstext}%
 \usepackage{amssymb}%
 \usepackage{showkeys}%
 \usepackage{cite}

\newcommand{\R}{\mathbb{R}}%
\DeclareMathOperator*\dom{dom}%
\DeclareMathOperator*\cl{cl}%
\DeclareMathOperator*\co{co}%
\DeclareMathOperator*\core{core}%
\DeclareMathOperator*\inte{int}%
\DeclareMathOperator*\sqri{sqri}%

\title{Closedness type regularity conditions for surjectivity results involving the sum of two maximal monotone operators
\thanks{Research partially supported by DFG (German Research Foundation), project WA 922/1-3.}}

\author{Radu Ioan Bo\c t
\thanks {Faculty of Mathematics, Chemnitz University of Technology,
D-09107 Chemnitz, Germany, e-mail:
radu.bot@mathematik.tu-chemnitz.de.} \and Sorin-Mihai Grad \thanks
{Faculty of Mathematics, Chemnitz University of Technology,
D-09107 Chemnitz, Germany, e-mail:
sorin-mihai.grad@mathematik.tu-chemnitz.de.}}
 \date{ }
\begin{document}
 \maketitle
\textbf{Abstract.} In this note we provide regularity conditions of closedness type which guarantee some surjectivity
results concerning the sum of two maximal monotone operators by using representative functions. The first regularity condition we give guarantees
the surjectivity of the monotone operator $S(\cdot + p)+T(\cdot)$, where $p\in X$ and $S$ and $T$ are maximal monotone operators
on the reflexive Banach space $X$. Then, this is used to obtain sufficient conditions for
the surjectivity of $S+T$ and for the situation when $0$ belongs to the range of $S+T$.
Several special cases are discussed, some of them delivering interesting byproducts.\\



\section{Introduction}

The recent developments in treating monotone operators by means of
convex analysis brought interesting results for many problems
involving monotone operators, like the maximal monotonicity of the
sum of two maximal monotone operators (cf. \cite{b, bgw}), but not
only (see \cite{b, SS}). Other interesting problems treated in
some recent articles like \cite{b, Z} deal with ranges of sums of
 maximal monotone operators, giving sufficient conditions that
ensure that 0 belongs to the range of a sum, respectively
that such a sum is surjective. Problems like these arise in different
applications in fields like inverse problems, Fenchel-Rockafellar
and Singer-Toland duality schemes, Clarke-Ekeland least action
principle (all in \cite{at}), variational inequalities (in
\cite{b, mt}), Schr\"odinger equations and others (cf.
\cite{abt}), while in papers like \cite{mo, mt} algorithms for
determining where does the sum of two maximal monotone operators
take the value 0 are given. Surjectivity issues regarding maximal monotone operators are discussed also in recent works such as \cite{ML, MAS, RM}.

In this paper we give, by using representative functions,
conditions that characterize the fact that, for the
maximal monotone operators $S$ and $T$ defined on the reflexive
Banach space $X$ and $p\in X$, the monotone operator $S(\cdot + p)+T(\cdot)$ is surjective.
From this we deduce characterizations of the surjectivity of $S+T$ and of the situation when $0$ lies in the range of $S+T$.
As main results, we introduce weak closedness type regularity conditions that guarantee the validity of the mentioned results.
An example to underline the fact these regularity conditions are indeed weaker than the interiority type ones considered in the literature is also provided. As special cases we consider situations where $T$ is the normal cone of a nonempty closed convex set,
respectively when $S$ and $T$ are subdifferentials of proper convex lower semicontinuous functions. In this way we rediscover several results from the literature, like the celebrated Rockafellar's surjectivity theorem and we moreover deliver weak regularity conditions for some results known so far only under stronger hypotheses involving generalized interiors.

\section{Preliminaries}

First we present some notions and results from convex analysis
that are necessary in order to make the paper self-contained. Let
the nontrivial separated locally convex topological space $X$ and its
continuous dual space $X^*$. The dual of $X^*$ is said to be the \textit{bidual} of $X$, being denoted by $X^{**}$. If $X$ is normed, it can be identified with a subspace of $X^{**}$, and we denote the canonical image in $X^{**}$ of the element $x\in X$ by $x$, too. By $\langle x^*, x\rangle$ we denote the value of the
linear continuous functional $x^*\in X^*$ at $x\in X$. Moreover, consider the \textit{duality product} $c:X\times X^*\rightarrow \R$,
$c(x, x^*)=\langle x^*, x\rangle$. Denote the
\textit{indicator} function of $U\subseteq X$ by $\delta_U$ and its \textit{support} function by $\sigma_U$.

For a function $f:X\rightarrow \overline\R = \R\cup \{ \pm\infty\}$, we denote its
\textit{domain} by $\dom f=\{x\in X: f(x)< +\infty\}$.
We call $f$ \textit{proper} if $f(x)>-\infty$ for all $x\in X$ and
$\dom f\neq\emptyset$. The \textit{conjugate} function of $f$ is
$f^*:X^*\rightarrow \overline\R$, $f^*(x^*)=\sup\big\{\langle
x^*, x\rangle - f(x) : x\in X\big\}$.
For $x\in X$ such that $f(x)\in \R$ we define the
\textit{(convex) subdifferential} of $f$ at $x$ by $\partial f (x)=\{x^*\in
X^*: f(y)-f(x)\geq \langle x^*, y-x\rangle\ \forall y\in X\}$. When $f(x)\notin \R$ we take by convention $\partial f(x)=\emptyset$.
The subdifferential of the indicator function of a set $U\subseteq X$ is said to be the \textit{normal cone} of that set.
Between a function and its
conjugate there is \textit{Young's inequality} $f^*(x^*)+f(x)\geq \langle x^*, x\rangle$ for all $x\in X$ and all $x^*\in X^*$,
fulfilled as equality by a pair $(x, x^*)\in X\times X^*$ if and
only if $x^*\in \partial f(x)$. Denote also by $\cl f:X\rightarrow \overline\R$ the
largest lower semicontinuous function everywhere less than or equal to $f$, i.e. the
\textit{lower semicontinuous hull} of $f$, and by $\co f:X\rightarrow \overline\R$ the
largest convex function everywhere less than or equal to $f$, i.e. the
\textit{convex hull} of $f$.

When $f, g:X\rightarrow
\overline\R$ are proper, we have the \textit{infimal convolution}
of $f$ and $g$ defined by $f\square g:X\rightarrow \overline\R$, $f\square g(a)=\inf
\{f(x)+g(a-x): x\in X\}$. It is said to be \textit{exact} at $y\in X$ when the
infimum at $a=y$ is attained, i.e. there exists $x\in X$ such that $f\square g (y)= f(x)+g(y-x)$.
When an infimum or a supremum is attained we
write $\min$, respectively $\max$ instead of $\inf$ and $\sup$.

The next result can be derived from the proofs of \cite[Proposition 2.2 and Theorem 3.1]{BW}.\\

{\bf Proposition 1.} {\it Consider on $X^*$ a locally convex topology giving $X$ as its dual space.
Let the proper, convex and lower
semicontinuous functions $f, g:X\rightarrow \overline \R$ satisfying $\dom f\cap \dom g\neq \emptyset$ and $p^*\in X^*$.
Then $f^*\square g^*$ is lower semicontinuous at $p^*$ and exact at $p^*$
if and only if
$$\inf_{x\in X} \big[f(x)+g(x) - \langle p^*, x\rangle\big]=\max_{x^*\in
X^*}\{-f^*(x^*)-g^*(p^*-x^*)\}.$$}

\textit{Remark 1.} The continuity of any of $f$ and $g$ at a point of $\dom f\cap \dom g$ yields the fulfillment of the equivalent statements from Proposition 1.\\

Let us recall some notions and results involving monotone
operators. Further $X$ is a Banach space equipped with the norm
$\|\cdot\|$, while the norm on $X^*$ is $\|\cdot\|_*$.

A multifunction $T:X \rightrightarrows X^*$ is called a
\textit{monotone operator} provided that for any $x, y \in X$ one
has $\langle y^*-x^*, y-x\rangle \geq 0$ whenever $x^*\in Tx$
 and $y^*\in Ty$. The \textit{domain} of $T$ is $D(T)=\{x\in X: Tx\neq\emptyset\}$, while its
\textit{range} is $R(T)=\cup \{Tx: x\in X\}$. $T$ is
called \textit{surjective} if $R(T)=X^*$. A monotone operator $T:X
\rightrightarrows X^*$ is called \textit{maximal} when its
\textit{graph} $G(T)=\{(x, x^*): x\in X, x^*\in Tx\}$ is not
properly included in the graph of any other monotone operator $T':
X \rightrightarrows X^*$. The subdifferential of a proper convex
lower semicontinuous function on $X$ is a typical example of a
maximal monotone operator.

To a maximal monotone operator $T:X \rightrightarrows X^*$ one can attach the
 \textit{Fitzpatrick function}
$$\varphi _T:X\times X^*\rightarrow \overline\R,\ \varphi
_T(x, x^*)=\sup\big\{\langle y^*, x\rangle + \langle x^*, y\rangle
- \langle y^*, y\rangle:y^*\in Ty\big\},$$ which is proper
convex and weak$\times$weak$^*$-lower semicontinuous, and the so-called {\it Fitzpatrick family} of \textit{representative functions}
$${\cal F}_T = \left\{f_T: X\times X^*\rightarrow\overline \R \bigg|
\begin{array}{l}
f_T \mbox{ is convex and strong$\times$strong lower semicontinuous},\\
c\leq f_T, (x, x^*)\in G(T) \Rightarrow f_T(x, x^*)=(x, x^*)
\end{array}\right\}.
$$
The largest element of ${\cal F}_T$ is $\psi_T=\cl_{\|\cdot\|\times \|\cdot\|_*} \co (c+\delta_{G(T)})$.
We also have $\varphi_T(x, x^*)=(c+\delta_{G(T)})^*(x^*, x)= \psi_T^*(x^*, x)$ for all $(x, x^*)\in X\times X^*$.
For $f_T\in {\cal F}_T$, denote by $\hat f_T:X\times X^*\rightarrow \overline\R$ the function defined as
$\hat f_T(x, x^*)=f_T(x, -x^*)$, $x\in X$, $x^*\in X^*$. Note that $\hat f_T$ is proper, convex and lower semicontinuous and
$\hat f_T(x, x^*)\geq -\langle x^*,
x\rangle$ and $\hat f_T^*(x^*, x)= f_T^*(x^*, -x)$ for all $x\in X$ and all $x^*\in X^*$.

If $f:X\rightarrow\overline\R$ is a
proper convex lower semicontinuous function, then the function $(x,
x^*)\mapsto f(x)+f^*(x^*)$ is a representative function of the
maximal monotone operator $\partial f:X \rightrightarrows X^*$ and
we call it the \textit{Fenchel representative function}. If $f$ is
moreover sublinear, the only representative function associated to
$\partial f$ is the Fenchel one, which coincides in this case with the Fitzpatrick function of $\partial f$.
Other maximal monotone operators having only one representative function, the Fenchel one, are the normal cones of nonempty closed convex sets.

The following statement underlines the close connections between
the maximal monotone operators and their representative functions.\\

\textbf{Proposition 2.} {\it Let $T:X \rightrightarrows X^*$ be a maximal monotone
operator.
\begin{enumerate}
\item[(i)] $\varphi_T$ is the smallest element of the family ${\cal F}_T$;

\item[(ii)] if $f_T\in {\cal F}_T$ one has $f_T^*(x^*, x)\geq \langle x^*, x\rangle$ for all $(x, x^*)\in X\times X^*$;

\item[(iii)] if $f_T\in {\cal F}_T$ and $(x, x^*)\in X\times X^*$, $(x, x^*)\in G(T)$ if and only if $f_T(x, x^*)=\langle x^*, x\rangle$ and this is further equivalent to $f_T^*(x^*, x)=\langle x^*, x\rangle$.
\end{enumerate}}

\section{Surjectivity results for the sum of two maximal monotone operators}

In this main section we deal with the
surjectivity results announced in the introduction. Let $X$ be further a reflexive Banach space and $S$ and
$T$ be two maximal monotone operators defined on $X$. The first main statement of this note follows, after an observation needed in its proof.\\

\textit{Remark 2.} Let $p\in X$ and $p^*\in X^*$. Then $p^*\in R(S(p + \cdot)+T(\cdot))$
if and only if $(p, p^*)\in G(S)-G(-T)$, where $G(-T)=\{(x, x^*)\in X\times X^*: (x, -x^*)\in G(T)\}$.\\

\textbf{Theorem 3.} {\it Let $p\in X$ and $p^*\in X^*$. The following statements are equivalent\\

$(i)$ \hspace{0.5cm} $p^*\in R(S(p + \cdot)+T(\cdot))$;

\begin{center}
$(ii)$ \begin{tabular}{r|l}
& $\forall f_S\in {\cal F}_S$ $\forall f_T\in {\cal F}_T$ one has $\dom f_S\cap (\dom \hat f_T + (p, p^*))\neq\emptyset$ and
the function
\\
&  $f_S^* \square \big(\hat f_T^* + \langle (p^*, p), (\cdot, \cdot)\rangle\big)$ is lower semicontinuous at $(p^*, p)$ and exact at $(p^*, p)$;
\end{tabular}
\end{center}

\begin{center}
$(iii)$ \begin{tabular}{r|l}
 & $\exists f_S\in {\cal F}_S$ $\exists f_T\in {\cal F}_T$  with $\dom f_S\cap (\dom \hat f_T + (p, p^*))\neq\emptyset$ such that the function\\
& $f_S^* \square \big(\hat f_T^* + \langle (p^*, p), (\cdot, \cdot)\rangle\big)$ is
lower semicontinuous at $(p^*, p)$ and exact at $(p^*, p)$.
\end{tabular}
\end{center}}

\textbf{Proof.} Note first that the assertion ``$(ii) \Rightarrow (iii)$'' is immediate and one also has
\begin{equation}\label{E3}\big(\hat f_T (\cdot-p, \cdot-p^*)\big)^* = \hat f_T^* + \langle p^*, \cdot\rangle
+ \langle \cdot, p\rangle.
 \end{equation}

``$(iii) \Rightarrow (i)$'' Proposition 1 yields the equivalence of $(iii)$ to
\begin{eqnarray}\label{E6}\big(f_S + \hat f_T (\cdot-p, \cdot-p^*)\big)^* (p^*, p) & = & \min_{u^*\in X^*, u\in X} \big[
f^*_S(p^*-u^*, p-u)\nonumber \\
& & + \hat f_T^*(u^*, u) + \langle p^*, u\rangle + \langle u^*, p\rangle\big].\end{eqnarray}
Denoting by $(\bar u^*, \bar u)\in X^*\times X$ the point where this minimum is attained, we obtain, via Proposition 2,
\begin{eqnarray}\label{E1}
\big(f_S+\hat f_T (\cdot-p, \cdot-p^*)\big)^* (p^*, p) &=& f^*_S(p^*-\bar u^*, p-\bar u) + \hat f_T^*(\bar u^*, \bar u) + \langle p^*, \bar u\rangle + \langle \bar u^*, p\rangle \nonumber\\
&\geq & \!\langle p^* \!-\! \bar u^*, p\!-\!\bar u\rangle \!-\! \langle  \bar u^*, \bar u\rangle \!+\!
\langle p^*, \bar u\rangle \!+\! \langle \bar u^*, p\rangle \!=\! \langle p^*, p\rangle.\end{eqnarray}
But Proposition 2 yields for every $x\in X$ and $x^*\in X^*$
$$
(f_S + \hat f_T (\cdot-p, \cdot-p^*)) (x, x^*) \geq  \langle x^*, x\rangle + \langle -(x^*-p^*), x-p\rangle =
\langle x^*, p\rangle + \langle p^*, x\rangle - \langle p^*, p\rangle,$$
thus $\langle p^*, p\rangle \geq \langle x^*, p\rangle + \langle p^*, x\rangle - (f_S + \hat f_T (\cdot-p, \cdot-p^*)) (x, x^*)$.
Consequently,
\begin{equation}\label{E4}
\big(f_S + \hat f_T (\cdot-p, \cdot-p^*)\big)^* (p^*, p)\leq \langle p^*, p\rangle.
\end{equation}
Together with \eqref{E1} this yields
$$
\big(f_S + \hat f_T (\cdot-p, \cdot-p^*)\big)^* (p^*, p)= \langle p^*, p\rangle,
$$
and consequently the inequalities invoked to obtain \eqref{E1} must be fulfilled as equalities.
Therefore
\begin{equation}\label{E2}
f_S^*(p^*-\bar u^*, p-\bar u) = \langle p^* - \bar u^*, p-\bar u\rangle\ \mbox{ and }\
\hat f_T^*(\bar u^*, \bar u) = \langle - \bar u^*, \bar u\rangle.
\end{equation}
Having these, Proposition 2 yields then
$p^*-\bar u^*\in S(p-\bar u)$ and $\bar u^*\in T(-\bar
u)$, followed by $p^*\in S(p-\bar u)+T(-\bar u)$, i.e. $p^*\in R(S(p + \cdot)+T(\cdot))$.

``$(i) \Rightarrow (ii)$'' Whenever $f_S\in {\cal F}_S$, $f_T\in {\cal F}_T$, $(i)$ yields, via Remark 2, $(p, p^*) \in \dom f_S - \dom \hat f_T$, i.e. $\dom f_S\cap (\dom \hat f_T + (p^*, p))\neq\emptyset$.

For every $f_S\in {\cal F}_S$, $f_T\in {\cal F}_T$, $u\in X$ and $u^*\in X^*$ we have
$f_S^*(p^*-u^*, p-u) + \hat f_T^* (u^*, u) + \langle (p^*, p), (u, u^*)\rangle$ $\geq \langle p^*-u^*, p-u\rangle - \langle u^*, u\rangle +
\langle p^*, u\rangle + \langle u^*, p\rangle = \langle p^*, p\rangle$, consequently,
$f_S^* \square \big(\hat f_T^* + \langle (p^*, p), (\cdot, \cdot)\rangle\big)$ $(p^*, p) \geq \langle p^*, p\rangle$ and, since the function in the right-hand side is strong$\times$strong continuous its value at $(p^*, p)$ must be also smaller than $\cl \big(f_S^* \square \big(\hat f_T^* + \langle (p^*, p), (\cdot, \cdot)\rangle\big)\big)(p^*, p)$. But from \cite[Theorem 7.6]{RB} we know, via \eqref{E3}, that
 one has $\cl \big(f_S^* \square \big(\hat f_T^* + \langle (p^*, p), (\cdot, \cdot)\rangle\big)\big) = (f_S + \hat f_T (-(p^*, p) + (\cdot, \cdot)))^*$ and since \eqref{E4} always holds, it follows that $\cl \big(f_S^* \square \big(\hat f_T^*$ $+ \langle (p^*, p), (\cdot, \cdot)\rangle\big)\big)(p^*, p)\leq \langle p^*, p\rangle$. Consequently,
\begin{equation}\label{E5}
f_S^* \square \big(\hat f_T^* + \langle (p^*, p), (\cdot, \cdot)\rangle\big)(p^*, p) \geq \cl \big(f_S^* \square \big(\hat f_T^* + \langle (p^*, p), (\cdot, \cdot)\rangle\big)\big)(p^*, p) = \langle p^*, p\rangle.
\end{equation}
Since $p^*\in R(S(p + \cdot)+T(\cdot))$, there exist $(\bar u^*, \bar u)\in X^*\times X$ fulfilling \eqref{E2}.
Then $f_S^*(p^*-\bar u^*, p-\bar u) + \hat f_T^*(\bar u^*, \bar u) + \langle (p^*, p), (\bar u, \bar u^*)\rangle
= \langle p^*, p\rangle$, i.e.
$$f_S^* \square \big(\hat f_T^* + \langle (p^*, p), (\cdot, \cdot)\rangle\big)(p^*, p) =
f_S^*(p^*-\bar u^*, p-\bar u) + \hat f_T^*(\bar u^*, \bar u) + \langle (p^*, p), (\bar u, \bar u^*)\rangle
= \langle p^*, p\rangle,$$ therefore the exactness of the infimal convolution in $(ii)$ is proven, while the lower semicontinuity
follows via \eqref{E5}.\hfill{$\Box$}\\

From Theorem 3 we obtain the following surjectivity result.\\

\textbf{Corollary 4.} {\it For $p\in X$, one has $R(S(p + \cdot)+T(\cdot))=X^*$ if and only if

\begin{center}
\begin{tabular}{|l}
$\forall p^*\in X^*$ $\forall f_S\in {\cal F}_S$ $\forall f_T\in {\cal F}_T$ one has $\dom f_S\cap (\dom \hat f_T + (p, p^*))\neq\emptyset$ and
the function
\\
  $f_S^* \square \big(\hat f_T^* + \langle (p^*, p), (\cdot, \cdot)\rangle\big)$ is lower semicontinuous at $(p^*, p)$ and exact at $(p^*, p)$,
\end{tabular}
\end{center}
and this is further equivalent to

\begin{center}
\begin{tabular}{|l}
 $\forall p^*\in X^*$ $\exists f_S\in {\cal F}_S$ $\exists f_T\in {\cal F}_T$ with $\dom f_S\cap (\dom \hat f_T + (p, p^*))\neq\emptyset$ such that the\\
function $f_S^* \square \big(\! \hat f_T^* \! +\!  \langle\!  (p^*, p), (\cdot, \cdot)\! \rangle\! \big)$ is
lower semicontinuous at $(p^*, p)$ and exact at $(p^*, p)$.\\
\end{tabular}
\end{center}}

Starting from Corollary 4 we are able to introduce a \textit{sufficient condition} for the surjectivity of $S(p + \cdot)+T(\cdot)$ for a given $p\in X$.\\

\textbf{Theorem 5.} {\it Let $p\in X$. Then $R(S(p + \cdot)+T(\cdot))=X^*$ if

\begin{center}
$(RC)$ \!\!\! \begin{tabular}{|l}\!\!\!\!
$\forall p^*\in X^*$ $\exists f_S\in {\cal F}_S$ $\exists f_T\in {\cal F}_T$ with $\dom f_S\cap (\dom \hat f_T + (p, p^*))\neq\emptyset$ such
that the func-\\
\!\!\!\! tion $f_S^* \square \big(\!\hat f_T^* \!+\! \langle (p^*, p), (\cdot, \cdot)\rangle\!\big)$ is
lower semicontinuous on $X^*\times \{p\}$ and exact at $(p^*, p)$.
\end{tabular}
\end{center}}

Next we characterize the surjectivity of $S+T$ via a condition involving representative functions.
The first statement follows directly from Theorem 3, while the second one is a direct consequence.\\

\textbf{Theorem 6.} {\it Let $p^*\in X^*$. The following statements are equivalent\\

$(i)$ \hspace{0.5cm} $p^*\in R(S+T)$;

\begin{center}
$(ii)$ \begin{tabular}{r|l}
& $\forall f_S\in {\cal F}_S$ $\forall f_T\in {\cal F}_T$ one has $\dom f_S\cap (\dom \hat f_T + (0, p^*))\neq\emptyset$ and
the function
\\
&  $f_S^* \square \big(\hat f_T^* + \langle p^*, \cdot\rangle\big)$ is lower semicontinuous at $(p^*, 0)$ and exact at $(p^*, 0)$;
\end{tabular}
\end{center}

\begin{center}
$(iii)$ \begin{tabular}{r|l}
 & $\exists f_S\in {\cal F}_S$ $\exists f_T\in {\cal F}_T$  with $\dom f_S\cap (\dom \hat f_T + (0, p^*))\neq\emptyset$ such that the function\\
& $f_S^* \square \big(\hat f_T^* + \langle p^*, \cdot\rangle\big)$ is
lower semicontinuous at $(p^*, 0)$ and exact at $(p^*, 0)$.\\
\end{tabular}
\end{center}}

\textbf{Corollary 7.} {\it One has $R(S+T)=X^*$ if and only if

\begin{center}
\begin{tabular}{|l}
$\forall p^*\in X^*$ $\forall f_S\in {\cal F}_S$ $\forall f_T\in {\cal F}_T$ one has $\dom f_S\cap (\dom \hat f_T + (0, p^*))\neq\emptyset$ and
the
\\
function $f_S^* \square \big(\hat f_T^* + \langle p^*, \cdot\rangle\big)$ is lower semicontinuous at $(p^*, 0)$ and exact at $(p^*, 0)$,
\end{tabular}
\end{center}
and this is further equivalent to

\begin{center}
\begin{tabular}{r|l}
 & $\forall p^*\in X^*$ $\exists f_S\in {\cal F}_S$ $\exists f_T\in {\cal F}_T$  with $\dom f_S\cap (\dom \hat f_T + (0, p^*))\neq\emptyset$ such that the \\
&   function $f_S^* \square \big(\hat f_T^* + \langle p^*, \cdot\rangle\big)$ is
lower semicontinuous at $(p^*, 0)$ and exact at $(p^*, 0)$.\\
\end{tabular}
\end{center}}

From Corollary 7 one can deduce a sufficient condition to have $S+T$ surjective.\\

\textbf{Theorem 8.} {\it One has $R(S+T)=X^*$ if

\begin{center}
$(\overline{RC})$ \!\!\! \begin{tabular}{|l}\!\!\!
 $\forall p^*\in X^*$ $\exists f_S\in {\cal F}_S$ $\exists f_T\in {\cal F}_T$  with $\dom f_S\cap (\dom \hat f_T + (0, p^*))\neq\emptyset$ such that the  \\
\!\!\! function $f_S^* \square \big(\hat f_T^* + \langle p^*, \cdot\rangle\big)$ is
lower semicontinuous on $X^*\times \{0\}$ and exact at $(p^*, 0)$.\\
\end{tabular}
\end{center}}

\textit{Remark 3.} In the literature there were given other regularity conditions guaranteeing the surjectivity
of $S+T$, namely, for fixed $f_S\in {\cal F}_S$ and $f_T\in {\cal F}_T$,
\begin{enumerate}
\item [-] (cf. \cite[Corollary 2.7]{ML}) $\dom f_T=X\times X^*$,
\item [-] (cf. \cite[Theorem 30.2]{SS}) $\dom f_S-\dom \hat f_T=X\times X^*$,
\item [-] (cf. \cite[Corollary 4]{Z}) $\{0\}\times X^*\subseteq \sqri(\dom f_S-\dom \hat f_T)$,
\end{enumerate}
where $\sqri$ denotes the \textit{strong quasi relative interior} of a given set, respectively. It is obvious that the first one implies the second, whose fulfillment yields the validity of the third condition. This one yields
$$\big(f_S + \hat f_T (\cdot, \cdot-p^*)\big)^* (x^*, 0) = \min_{u^*\in X^*, u\in X} \big[
f^*_S(x^*-u^*, -u)+ \hat f_T^*(u^*, u) + \langle p^*, u\rangle \big]\ \forall x^*, p^*\in X^*,$$
which is equivalent, when $\dom f_S\cap (\dom \hat f_T + (0, p^*))\neq\emptyset$ (condition fulfilled in all the three regularity conditions given above), to the fact that whenever $p^*\in X^*$ the function
$f_S^* \square \big(\hat f_T^* + \langle p^*, \cdot\rangle\big)$ is lower semicontinuous at $(x^*, 0)$ and exact at $(x^*, 0)$ for all $x^*\in X^*$. It is obvious that this implies $(\overline{RC})$ and below we present a situation where $(\overline{RC})$ holds, unlike the conditions cited from the literature for the surjectivity of $S+T$.\\

\textit{Example 1.} Let $X=\R$ and consider the operators
$S, T: \R\rightrightarrows \R$ defined by
$$Sx=\left\{\begin{array}{ll}
\{0\},& \mbox{if } x>0,\\
(-\infty, 0],& \mbox{if } x=0,\\
\emptyset,& \mbox{otherwise},\end{array}\right. \ \mbox{ and }
Tx=\left\{\begin{array}{ll} \R,& \mbox{if } x=0,\\ \emptyset,&
\mbox{otherwise},\end{array}\right. x\in \R.$$
One notices
easily that, considering the functions $f,
g:\R\rightarrow\overline\R$, $f=\delta_{[0, +\infty)}$ and
$g=\delta_{\{0\}}$, which are proper, convex and
lower-semicontinuous, we have $S=\partial f$ and $T=\partial g$, thus $S$ and $T$ are maximal monotone.
It is obvious that $R(S+T)=\R$.
The Fitzpatrick families of both $S$ and $T$ contain only their Fitzpatrick functions, because $f$ and $g$ are sublinear functions.
We have
$$\varphi_S(x, x^*)=\left\{\begin{array}{ll}
0,& \mbox{if } x\geq 0,  x^*\leq 0,\\
+\infty,& \mbox{otherwise},\end{array}\right. \mbox{and }
\varphi_T(x, x^*)=\left\{\begin{array}{ll}
0,& \mbox{if } x=0,\\
+\infty,& \mbox{otherwise},\end{array}\right.$$
therefore
$$\varphi_S^*(x^*, x)=\left\{\begin{array}{ll}
0,& \mbox{if } x^*\leq 0, x\geq 0,\\
+\infty,& \mbox{otherwise},\end{array}\right. \mbox{and }
\varphi_T^*(x^*, x)=\left\{\begin{array}{ll}
0,& \mbox{if } x=0,\\
+\infty,& \mbox{otherwise}.\end{array}\right.$$
Then $\dom \varphi_S-\dom \hat\varphi_T = \R_+\times \R$, where $\R_+=[0, +\infty)$, and it is obvious that
$\{0\}\times \R$ is not included in $\sqri(\dom \varphi_S-\dom \hat \varphi_T) = (0, +\infty)\times\R$.
Consequently, the three conditions mentioned in Remark 3 fail in this situation. On the other hand, for $p^*, x, x^*\in \R$ one has
$$\varphi_S^*\square \big(\hat \varphi_T^* + \langle p^*, \cdot\rangle\big) (x^*, x) =
\left\{\begin{array}{ll}
0,&\mbox{if } x\geq 0,\\
+\infty,& \mbox{ if } x< 0.
\end{array}\right.$$
This function is lower semicontinuous on $\R\times \R_+$ and exact at all $(x^*, x)\in \R\times \R_+$. Consequently, $(\overline{RC})$ is valid in this case.\\

\textit{Remark 4.} Following Remark 1, when one of $f_S$ and $f_T$ is continuous, $(RC)$ is automatically fulfilled.
It is known (see for instance \cite{SS}) that the domain of the Fitzpatrick function attached to the \textit{duality map}
$$J:X\rightrightarrows X^*,\ Jx=\partial\frac
{1}{2}\|x\|^2=\Big\{x^*\in X^*: \|x\|^2 = \|x^*\|_*^2 = \langle x^*,
x\rangle\Big\}, x\in X,$$ which is a maximal monotone
operator, is the whole product space $X\times X^*$. By \cite[Theorem 2.2.20]{ZZ} it follows that $\varphi_J$ is continuous,
thus by Corollary 4 we obtain that $S(p+\cdot) + J(\cdot)$ is surjective, whenever $p\in X$.
Thus we rediscover a property of the maximal monotone operators. Moreover, employing Corollary 7 one gets that $S+J$ is surjective,
rediscovering Rockafellar's classical surjectivity theorem for maximal monotone operators (cf. \cite[Theorem 29.5]{SS}).\\

The last results we derive from the main one are connected to the situation when 0 lies in the range of $S+T$.\\

\textbf{Corollary 9.} {\it One has $0\in R(S+T)$ if and only if
\begin{center}
\begin{tabular}{|l}
 $\forall f_S\in {\cal F}_S$ $\forall f_T\in {\cal F}_T$ one has $\dom f_S\cap \dom \hat f_T \neq\emptyset$
and the function\\
  $f_S^* \square \hat f_T^*$ is
lower semicontinuous at $(0, 0)$ and exact at $(0, 0)$,
\end{tabular}
\end{center}
and this is further equivalent to
\begin{center}
\begin{tabular}{r|l}
 & $\exists f_S\in {\cal F}_S$ $\exists f_T\in {\cal F}_T$  with $\dom f_S\cap \dom \hat f_T \neq\emptyset$ such that the function\\
&   $f_S^* \square \hat f_T^*$ is
lower semicontinuous at $(0, 0)$ and exact at $(0, 0)$.
\end{tabular}
\end{center}}

From Corollary 9 one can deduce a sufficient condition to be sure that $0\in R(S+T)$.\\

\textbf{Theorem 10.} {\it One has $0\in R(S+T)$ if
\begin{center}
$(\widetilde{RC})$
\begin{tabular}{|l}
 $\exists f_S\in {\cal F}_S$ $\exists f_T\in {\cal F}_T$  with $\dom f_S\cap \dom \hat f_T \neq\emptyset$ such that the function\\
  $f_S^* \square \hat f_T^*$ is
lower semicontinuous on $X^*\times \{0\}$ and exact at $(0, 0)$.\\
\end{tabular}
\end{center}}

\textit{Remark 5.} Other regularity conditions guaranteeing $0\in R(S+T)$
were given in \cite[Theorem 4.5]{b}, $(0, 0) \in \core (\co(G(S))-\co(G(-T)))$, where
$\core$ denotes the \textit{algebraic interior} of a given set and $\co$ its \textit{convex hull},
and \cite[Lemma 1]{Z}, $(0, 0)\in \sqri(\dom f_S-\dom \hat f_T)$, respectively. Following similar arguments to the ones in Remark 3 one can
show that both yield the validity of $(\widetilde{RC})$. Checking the situation from Example 1, we see that the condition involving $\sqri$ fails,
while $(\widetilde{RC})$ is fulfilled. As $\core (\co(G(S))-\co(G(-T)))= \inte (\R_+\times (-\R_+) - \{0\}\times \R) = (0, +\infty)\times \R$ does not contain $(0, 0)$, it is straightforward that $(\widetilde{RC})$ is indeed weaker than both the mentioned conditions for $0\in R(S+T)$ from the literature.\\

\textit{Remark 6.} One can notice via \eqref{E3} that \eqref{E6} can be rewritten when $p^*=0$ and $p=0$ as
\begin{equation}\label{E7}
\inf_{x\in X, x^*\in X^*}
\big[f_S (x, x^*) + \hat f_T  (x, x^*)\big] = \max_{u^*\in X^*, u\in X} \big[- f^*_S(-u^*, -u) - \hat f_T^*(u^*, u) \big],
\end{equation}
i.e. there is strong duality for the convex optimization problem
formulated above in the left-hand side of \eqref{E7} and its
Fenchel dual problem. When $(\bar u, \bar u^*)\in X\times X^*$ is an optimal solution of the dual problem, i.e.
the point where the maximum in the right-hand side of \eqref{E7} is attained, one obtains
$\bar u^*\in S(\bar u)$ and $-\bar u^*\in T(\bar u)$. Employing now Proposition 2 we obtain
$f_S (\bar u, \bar u^*)=f_S^* (-\bar u^*, -\bar u)=\langle \bar u^*, \bar u\rangle$ and $\hat f_T (\bar u, \bar u^*)=
\hat f_T^* (\bar u^*, \bar u)=-\langle \bar u^*, \bar u\rangle$, therefore
$$f_S (\bar u, \bar u^*) + \hat f_T (\bar u, \bar u^*) = f_S^* (-\bar u^*, -\bar u) + \hat f_T^* (\bar u^*, \bar u)=
0.$$
Thus, the infimum in the left-hand side of \eqref{E7} is attained, i.e. the primal optimization problem
given there has an optimal solution, too. As strong duality for it holds, we are now in the situation called \textit{total duality} (cf. \cite[Section 17]{RB}), which happens when the optimal objective values of the primal and dual coincide and both these problems have optimal solutions.
Therefore we can conclude that for this kind of optimization problems when strong Fenchel duality holds the primal problem has an optimal solution, too.\\

\textit{Remark 7.} Given $p\in X$ and $p^*\in X^*$, the function
$f_S^* \square \big(\hat f_T^* + \langle (p^*, p), (\cdot, \cdot)\rangle\big)$ can be replaced in
conditions $(ii)-(iii)$ from Theorem 3 with
$\big(f_S^* - \langle (p^*, p), (\cdot, \cdot)\rangle\big) \square \hat f_T^*$ without altering the statement.
The other conditions considered above can be correspondingly rewritten, too.

\section{Applications}

\subsection{When $T$ is the normal cone of a closed convex set}

Let $U\subseteq X$ be a nonempty closed convex set. Its normal
cone $N_U$ is a maximal monotone operator.
Taking $T=N_U$, its only representative function is $f_{N_U}(x,
x^*)=\delta_U(x)+\sigma_U(x^*)$, $(x, x^*)\in X\times X^*$. From our main statements we obtain in this case the following results.\\

\textbf{Corollary 11.} {\it Let $p\in X$. Then $R(S(p + \cdot)+N_U(\cdot))=X^*$ if and only if
\begin{center}
\begin{tabular}{r|l}
& $\forall p^*\in X^*$ $\forall f_S\in {\cal F}_S$ one has $\dom f_S\cap (U\times \dom \sigma_{-U} + (p, p^*))\neq\emptyset$ and
the function
\\
&  $(y^*, y)\mapsto \inf_{x\in -U, x^*\in X^*}\big[(f_S^* -  \langle (p^*, p), (\cdot, \cdot)\rangle)(y^*-x^*, y-x) + \sigma_U(x^*)\big]$ is lower\\
 &  semicontinuous at $(p^*, p)$ and the infimum within is attained when $(y^*, y)=(p^*, p)$,
\end{tabular}
\end{center}
and this is further equivalent to
\begin{center}
\begin{tabular}{r|l}
& $\forall p^*\in X^*$ $\exists f_S\in {\cal F}_S$ with $\dom f_S\cap (U\times \dom \sigma_{-U} + (p, p^*))\neq\emptyset$
the function
\\
&  $(y^*, y)\mapsto \inf_{x\in -U, x^*\in X^*}\big[(f_S^* -  \langle (p^*, p), (\cdot, \cdot)\rangle)(y^*-x^*, y-x) + \sigma_U(x^*)\big]$ is lower\\
 & semicontinuous at $(p^*, p)$ and the infimum within is attained when $(y^*, y)=(p^*, p)$.\\
\end{tabular}
\end{center}}

\textbf{Corollary 12.} {\it Let $p\in X$. Then $R(S(p + \cdot)+N_U(\cdot))=X^*$ if
\begin{center}
\begin{tabular}{r|l}
&\!\! $\forall p^*\in X^*$ $\exists f_S\in {\cal F}_S$ with $\dom f_S\cap (U\times \dom \sigma_{-U} + (p, p^*))\neq\emptyset$
the function\\
$(RC_U)$\!\!\! &\!\! $(y^*, y)\mapsto \inf_{x\in -U, x^*\in X^*}\big[(f_S^* -  \langle (p^*, p), (\cdot, \cdot)\rangle)(y^*-x^*, y-x) + \sigma_U(x^*)\big]$ is lower\\
&\!\! semicontinuous on $X^*\times\{p\}$ and the infimum within is attained when $(y^*, y)=(p^*, p)$.\\
\end{tabular}
\end{center}}

\textbf{Corollary 13.} {\it One has $0\in R(S+N_U)$ if
\begin{center}
$(RC_0)$
\begin{tabular}{|l}
 $\exists f_S\in {\cal F}_S$ with $\dom f_S\cap (U\times \dom \sigma_{-U} )\neq\emptyset$ such that the function\\
   $(y^*, y)\mapsto \inf_{x\in U} \{(f_S^*(\cdot, y+x)\square \sigma_U)(y^*)\}$ is lower semicontinuous \\
 on $X^*\times\{0\}$ and the infimum within is attained when $(y^*, y)=(0, 0)$.\\
\end{tabular}
\end{center}}

\textit{Remark 8.} A stronger than $(RC_0)$ regularity condition for $0\in R(S+N_U)$, namely $0\in \core (\co(D(S)$ $-U))$,
was considered in  \cite[Corollary 5.7]{b}.\\

Not without importance is the question how can one equivalently characterize the surjectivity of a maximal monotone operator via its representative
functions. To proceed to answering it, take $U=X$. Then $T=N_X$, i.e. $Tx=\{0\}$ for all $x\in X$, and the Fenchel representative function of $N_X$ is
$(x, x^*)\mapsto \delta_X(x) + \sigma_X(x^*) = \delta_{0}(x^*)$. Then $S+T=S$ and the surjectivity of $S$ can be characterized, via Corollary 11, as follows.\\

\textbf{Corollary 14.} {\it One has $R(S)=X^*$ if and only if
\begin{center}
\begin{tabular}{r|l}
& $\forall p^*\in X^*$ $\forall f_S\in {\cal F}_S$ the function $y^*\mapsto -\big(f_S^* (y^*, \cdot)\big)^*(p^*)$
is lower \\
&  semicontinuous at $p^*$ and $\exists x\in X$ such that $p^*\in (\partial f_S^* (p^*, \cdot))(x)$,
\end{tabular}
\end{center}
and this is further equivalent to
\begin{center}
\begin{tabular}{r|l}
& $\forall p^*\in X^*$ $\exists f_S\in {\cal F}_S$ the function $y^*\mapsto -\big(f_S^* (y^*, \cdot)\big)^*(p^*)$
is lower \\
&  semicontinuous at $p^*$ and $\exists x\in X$ such that $p^*\in (\partial f_S^* (p^*, \cdot))(x)$.
\end{tabular}
\end{center}}

\textbf{Proof.} Corollary 11 asserts the equivalence of the surjectivity of the maximal monotone operator $S$ to the lower
semicontinuity at $(p^*, 0)$ of the function $(y^*, y)\mapsto \inf_{x\in X, x^*\in X^*}\big[(f_S^* -  \langle p^*, \cdot\rangle)(y^*-x^*, y+x) + \sigma_X(x^*)\big]$ concurring with the attainment of the infimum within when $(y^*, y)=(p^*, 0)$ for every $p^*\in X^*$.

Taking a closer look at this function, we note that it can be simplified to
$(y^*, y)\mapsto \inf_{x\in X}$ $\big[f_S^* (y^*, y+x) - \langle p^*,  y+x\rangle\big]$, which can be further reduced to
$y^*\mapsto -\big(f_S^* (y^*, \cdot)\big)^*(p^*)$.

For $p^*\in X^*$, the attainment of the infimum from above when $(y^*, y)=(p^*, 0)$ means actually the existence of $x\in X$ such that $f_S^* (p^*, x) - \langle p^*,  x\rangle = -\big(f_S^* (p^*, \cdot)\big)^*(p^*)$, which is nothing but $p^*\in (\partial f_S^* (p^*, \cdot))(x)$.\hfill{$\Box$}\\

\textit{Remark 9.} In \cite[Corollary 2.2]{ML} it is said that $S$ is surjective if $\dom(\varphi_S)=X\times X^*$. This result can be obtained as a
consequence of Corollary 14, via Remark 1.\\

\textit{Remark 10.} Determining when $0\in R(S)$ is important even beyond optimization.
Via Corollary 14 we can provide the following sufficient condition for this
\begin{center}
\begin{tabular}{r|l}
& $\exists f_S\in {\cal F}_S$ the function $y^*\mapsto -\big(f_S^* (y^*, \cdot)\big)^*(0)$
is lower  \\
& semicontinuous and $\exists x\in X$ such that $p^*\in (\partial f_S^* (0, \cdot))(x)$.
\end{tabular}
\end{center}

\subsection{When $S$ and $T$ are subdifferentials}

Let the proper convex lower semicontinuous functions $f, g:X\rightarrow \overline \R$. Take first $T=\partial g$ and consider for it the Fenchel representative function. Then Corollary 4 yields the following statement.\\

\textbf{Corollary 15.} {\it Let $p\in X$. Then $R(S(p+\cdot)+\partial g(\cdot))=X^*$ if and only if
\begin{center}
\begin{tabular}{r|l}
& $\forall p^*\in X^*$ $\forall f_S\in {\cal F}_S$ one has $\dom f_S\cap (\dom g \times (-\dom g^*) + (p, p^*))\neq\emptyset$ and
the function
\\
&  $f_S^* \square \big(g(-\cdot) + g^*(\cdot)  + \langle (p^*, p), (\cdot, \cdot)\rangle\big)$ is lower semicontinuous at $(p^*, p)$ and exact at $(p^*, p)$,\\
\end{tabular}\\
\end{center}

and this is further equivalent to
\begin{center}
\begin{tabular}{r|l}
 & $\forall p^*\in X^*$ $\exists f_S\in {\cal F}_S$ with $\dom f_S\cap (\dom g \times (-\dom g^*) + (p, p^*))\neq\emptyset$ such that the function\\
&  $f_S^* \square \big(g(-\cdot) + g^*(\cdot)  + \langle (p^*, p), (\cdot, \cdot)\rangle\big)$ is lower semicontinuous at $(p^*, p)$ and exact at $(p^*, p)$.\\
\end{tabular}
\end{center}}

The other statements in Section 3 can be particularized for this special case, too. However, we give here only a consequence of Theorem 10.\\

\textbf{Corollary 16.} {\it One has $0\in R(S+\partial g)$ if
\begin{center}
\begin{tabular}{r|l}
 & $\exists f_S\in {\cal F}_S$ with $\dom f_S\cap (\dom g \times (-\dom g^*))\neq\emptyset$ such that the function\\
&  $f_S^* \square \big(g(-\cdot) + g^*(\cdot)\big)$ is lower semicontinuous  on $X^*\times\{0\}$ and exact at $(0, 0)$.\\
\end{tabular}
\end{center}}

\textit{Remark 11.} In \cite[Proposition 2.9]{ML} it was proven that when $g$ and $g^*$ are real-valued the monotone operator $S(p+\cdot)+\partial g(\cdot)$ is surjective whenever $p\in X$. This statement can be rediscovered as a consequence of Corollary 15, too. Using \cite[Proposition 2.1.6]{ZZ} one obtains that $g$ and $g^*$ are continuous under the mentioned hypotheses. Then the Fenchel representative function of $\partial g$ is continuous and (see Remark 1) this yields the fulfillment of the regularity condition from Corollary 15. Consequently, $S(p+\cdot)+\partial g(\cdot)$ is surjective whenever $p\in X$.\\

Take now also $S=\partial f$, to which we associate the Fenchel representative function, too. Let the function $\hat g:X\rightarrow \overline\R$,  $\hat g(x)=g(-x)$. Corollary 4 yields the following result.\\

\textbf{Corollary 17.} {\it Let $p\in X$. If $\dom f\cap (p+\dom g)\neq \emptyset$, then $R(\partial f(p+\cdot)+\partial g(\cdot))=X^*$ if and only if
\begin{center}
\begin{tabular}{|l}
$\forall p^*\in X^*$ one has $\dom f^*\cap (p^*-\dom g^*)\neq\emptyset$, the function $f\square (\hat g + p^*)$ \\
 is lower semicontinuous at $p$ and exact at $p$ and the function\\
$f^* \square (g^* + p)$ is lower semicontinuous at $p^*$ and exact at $p^*$.\\
\end{tabular}
\end{center}}

Moreover, from Corollary 16 one can deduce the following statement.\\

\textbf{Corollary 18.} {\it One has $0\in R(\partial f+\partial g)$ if $\dom f\cap \dom g\neq \emptyset$, $\dom f^*\cap (-\dom g^*) \neq\emptyset$ and
\begin{center}
\begin{tabular}{r|l}
& $f\square \hat g $ is lower semicontinuous at $0$ and exact at $0$ and the \\
& function $f^* \square g^*$ is lower semicontinuous and exact at $0$.
\end{tabular}
\end{center}}

\end{document}